\documentclass[11pt]{article}
\usepackage{mathrsfs}
\usepackage{amssymb}
\usepackage{amsfonts}
\usepackage{color,xcolor}
\usepackage{graphicx}
\usepackage{geometry}
\usepackage{manfnt}
\usepackage[pagewise]{lineno}
\RequirePackage[colorlinks,citecolor=blue,urlcolor=blue]{hyperref}
\newtheorem{theorem}{Theorem}[section]
\newtheorem{corollary}{Corollary}[section]
\newtheorem{lemma}{Lemma}[section]

\newtheorem{definition}{Definition}[section]
\newtheorem{remark}{Remark}[section]

\newcommand{\beq}{\begin{equation}}
\newcommand{\eeq}{\end{equation}}
\newcommand{\beqn}{\begin{eqnarray}}
\newcommand{\eeqn}{\end{eqnarray}}

\begin{document}
\title{\bf The central configuration of the planar ($N$+1)-body problem with a regular $N$-polygon}
\author{Liang Ding$^{1}$, Jinlong Wei$^{2}$ and Shiqing Zhang$^{3\,\ast}$
\\ {\small \it $^1$School of Data Science and Information Engineering, Guizhou} \\  {\small \it Minzu University, Guiyang, 550025, China}\\
 {\small \tt ding2016liang@126.com}\\ {\small \it $^2$School of Statistics and Mathematics, Zhongnan University} \\  {\small \it of Economics and Law, Wuhan, 430073, China}\\
 {\small \tt weijinlong.hust@gmail.com}\\ {\small \it $^3$School of Mathematics, Sichuan University} \\  {\small \it  Chengdu, 610064, China}\\
 {\small \tt zhangshiqing@msn.com}}
\date{}
 \maketitle
\noindent{\hrulefill} \vskip1mm\noindent
 {\bf Abstract}
For planar ($N$+1)-body ($N$\,$\geq$ 2) problem with a regular $N$-polygon, under the assumption that the ($N$+1)-th body locates at the geometric center of the regular $N$-polygon, we obtain the sufficient and necessary conditions that the $N$+1 bodies can form a central configuration.



 \vskip0mm\noindent{\hrulefill}

\section{Introduction} \label{sec1}
\setcounter{equation}{0}
We have known that the motion of Newtonian ($N$+1)-body problem with positive masses $m_{k}$ and positions $x_{k}(t)$ is described by Newton's laws of motion:
\begin{eqnarray}\label{1.1}
m_k\ddot{x}_k(t)=\sum_{{j\neq k\atop{{1\leq j\leq N+1}}}}m_km_j\frac{x_j(t)-x_k(t)}{|x_j(t)-x_k(t)|^3}, \quad x_{k}(t)\in\mathbb{R}^{d}, \quad d=2 \,\,or\,\,3, \,\,k=1,2,\ldots,N+1,
\end{eqnarray}
and when $d=2$, we call it planar ($N$+1)-body problem. Since finding the planar central configurations and the relative equilibrium solutions of the system (\ref{1.1}) are equivalent \cite{Saari2005}, then before introducing the definition on central configuration in $\mathbb{R}^{2}$, we introduce the definition of relative equilibrium solution firstly.
\begin{definition}\cite[Definition 2.3.2]{Llibre2015}
For a solution $\vec{x}(t)=(x_{1}(t),x_{2}(t),...,x_{N+1}(t))^{T}$ of system (\ref{1.1}), we call it
a relative equilibrium solution if
$x_{j}(t)$ has the form that $x_{j}(t)=(x_{j}(0)-c_{0})e^{i\omega t}+c(t)$ where $x_{j}(0)$,
$c_{0}=\sum\limits_{ j=1}^{N+1}m_jx_{j}(0)/\sum\limits_{j=1}^{N+1}m_{j}$,
$c(t)=\sum\limits_{j=1}^{N+1}m_jx_{j}(t)/\sum\limits_{j=1}^{N+1}m_{j}$ and $\omega$, represent the initial position of the $j$-th body ($j=1,2,\ldots,N+1$), the center of masses of the $N$+1 bodies at initial instant, the center of masses of the $N$+1 bodies at any instant, and a uniform angular velocity of rotation around $c_{0}$, respectively.
\end{definition}

It is well-known that the above relative equilibrium solution $\vec{x}(t)$ exists if and only if
$x_{j}(0)\in\mathbb{R}^{2}$ satisfy the following algebraic equations
\begin{eqnarray}\label{1.2}
\sum\limits_{{j\neq k\atop{{1\leq j\leq N+1}}}}
\frac{m_j m_k}{|x_{j}(0)-x_{k}(0)|^{3}}(x_{j}(0)-x_{k}(0))=-\omega^{2} m_{k}(x_{k}(0)-c_{0}), \quad \  \  k=1,2,...,N+1,
\end{eqnarray}
Employing (\ref{1.2}), we can give the definition of planar central configuration.
\begin{definition}
A planar configuration $\vec{x}(0)=\{(x_{1}(0),x_{2}(0),...,x_{N+1}(0))^{T}\in (\mathbb{R}^{2})^{N+1}: x_{k}(0)\neq x_{j}(0) \,\, when \,\, k\neq j\}$, is called a planar central configuration if $\vec{x}(0)$ satisfies (\ref{1.2}).
\end{definition}
In fact, in (\ref{1.2}), if we substitute $x_{j}(0)\in\mathbb{R}^{3}$ and general constant $\lambda$ for $x_{j}(0)\in\mathbb{R}^{2}$ and $-\omega^{2}$ respectively, then we call $\vec{x}(0)$ a general central configuration, and the study of central configurations is a very important and difficult topic in celestial mechanics with a long and varied history \cite{Wintner1947}, especially, the problem on the numbers for central configurations is so important that, in 2000, Smale \cite{Smale2000} took it as one of the most important 18 mathematical problems for the 21st century. Here we only mention some related
works to this paper. For 3-body problem ($N+1=3$), in 1772, Lagrange \cite{Lagrange1772} established the well known `equilateral-triangle solutions'.
For $N+1\geq4$, in 1985, Perko and Walter \cite{Perko1985} proved that if
$\vec{x}(0)=(x_1(0),...,x_{N+1}(0))^{T}$ located at the vertices of planar regular polygon, then they formed a regular polygonal central configuration if and only if all the masses were equal. Moreover, in 2010, by using analytic methods, Shi and Xie \cite{Shi2010} arrived at the conclusion that for the central configuration of planar four body problem consisting of three particles of equal mass, besides the family of equilateral triangle configurations, there are exactly one family of concave and one family of convex central configurations. In 2013, Li and Wang \cite{Li2013} considered central configurations of $N$-body problem, and they showed that if $N$ point masses were located at the vertices of a regular ($N$+1)-polygon, then the $N$ particles did not form a central configuration for any value of masses. In 2018, for the planar $N$+1-body problem, Fernandes, Garcia, Llibre and Mello \cite{Fernandes2018} studied the central configuration that $N$ equal masses locates at the vertices of a regular polygonal and the ($N$+1)-th body with null mass at the plane, and they obtained that when $N\geq3$, there are $3N$+1 classes of central configurations. For other related works on regular polygonal central configurations, we refer to \cite{Chen2018,Fay1996,Moeckel1995,Ouyang2004,Saari2005,Yu2012,Zhang-Zhou,Zhang-Zhu}, especially for spatial twisted regular polygonal central configurations, one consults to \cite{Moeckel1995,Yu2012,Zhang-Zhou,Zhang-Zhu}, and for pyramid central configurations, please see \cite{Fay1996,Ouyang2004}.

Note that finding the relative equilibrium solutions of the classical $N$-body problem and the planar central configurations are equivalent, and in the definition of relative equilibrium solution, $x_{j}(t)$ ($j=1,2,\ldots,N+1$) rotates around the mass center $c_{0}$, and the ($N$+1)-th body is moving, and we also note that in 2018, for the solution of the planar system (\ref{1.1}), under the assumption that the ($N$+1)-th ($N\geq4$) body is motionless at any instant, the other $N$ particles rotate around the ($N$+1)-th particle with the same angular velocity $\omega$, Chen and Luo \cite{Chen2018} proved that the $N$ bodies must have equal mass and ($N$+1)-th body must be at the mass center of the $N$+1 bodies. We notice that, in the general case, the mass center of the $N$+1 bodies may not coincide with the positions of the ($N$+1)-th body or the geometry center of all the $N$+1 bodies, then it is natural to ask that

\textbf{Question:} For the planar ($N$+1)-body problem, if the ($N$+1)-th body is moving, and the other $N$ particles rotate around the mass center $c_{0}$ of the $N$+1 bodies with the same angular velocity $\omega$, does the existence of the solution of system (\ref{1.1}) still implies that all the masses in the regular $N$-gon must be equal to each other ?

 In fact, by the relationship between relative equilibrium solution and planar central configuration, we see that the assumption of the above question is equivalent to that all the
$N$+1 bodies form a planar central configuration, then in this paper, we attempt to answer the above question.

\section{Main result} \label{sec2}

In the following, we assume that the given $N$+1 particles are in the same plane at the initial instant $t=0$. For $j=1,2,\ldots,N$, the particle $x_j(0)$ with positive mass $m_j$ locate at the vertices of the regular
$N$-polygon, and the
($N$+1)-th particle with positive mass $m_{N+1}$ locates at the geometric center of the regular $N$-polygon. Without loss of generality, we assume the coordinates of the particle $x_j(0)$ be $e^{i\theta_{j}}=e^{i\frac{2j\pi}{N}}=:q_j(0)$ for $j=1,2,\ldots,N$, where $q_{1}(0), q_{2}(0),\ldots,q_{N}(0)$ are the $N$ roots of unity. Then the main result is described as follows.

\begin{theorem}\label{thm2.1}
For the configuration of planar ($N$+1)-body ($N\geq2$) problem with a regular $N$-polygon, under the assumption that the ($N$+1)-th body locates at the geometric center of the regular $N$-polygon, we prove that the $N$+1 bodies can form a central configuration if and only if
\begin{eqnarray*}
m_{1}=m_{2}=\ldots=m_{N}=\frac{\omega^{2}-m_{N+1}}{\frac{1}{4}\sum_{j=1}^{N-1}\csc\frac{j\pi}{N}}.
\end{eqnarray*}
\end{theorem}

\begin{remark}\label{rem2.1}
In \cite{Chen2018},  under the assumptions that the $N$ particles located at the vertices of the regular
$N$-polygon, rotated around the ($N$+1)-th particle with the same angular velocity $\omega$, and the ($N$+1)-th ($N\geq4$) body is assumed motionless at any instant, Chen and Luo obtained the same conclusion as in Theorem \ref{thm2.1}.
But now, by using the equivalence between relative equilibrium solutions and planar central configuration, under the assumptions that
all the $N$+1 bodies rotate around the mass center $c_{0}$ of the $N$+1 bodies with the same angular velocity $\omega$, and the ($N$+1)-th body locates at the geometric center of the regular $N$-polygon, we obtain the same sufficient and necessary conditions for the $N$+1 bodies which can form a central configuration as \cite{Chen2018}. Moreover, notice that we can not know the geometric center coincides with the center of masses in advance, so we do not assume that the ($N$+1)-th ($N\geq2$) body is motionless at any instant.
\end{remark}

\begin{remark}\label{rem2.2}
For the 3-body problem ($N$+1=3), the well-known Eulerian collinear central configuration exists for any choice of the masses of three bodies. Therefore the assumption of Theorem \ref{thm2.1}, i.e  the ($N$+1)-th body locates at the geometric center of the regular $N$-polygon is necessary.
\end{remark}

Based upon Theorem \ref{thm2.1}, we have
\begin{corollary}\label{cor2.1}
For the central configuration of planar ($N$+1)-body ($N\geq2$) problem with a regular $N$-polygon, if the ($N$+1)-th body locates at the geometric center of the regular $N$-polygon, then the geometric center of the $N$+1 bodies is just the center of the masses of the $N$+1 bodies.
\end{corollary}

\section{Some Lemmas}\label{sec3}\setcounter{equation}{0}
At beginning, we introduce a result about Eulerian collinear central configuration:
\begin{lemma}\label{lem3.1}\cite[Page 94, lines 18-19 and Page 95, lines 1-2]{Siegel1971}
For 3-body problem ($N$+1=3), set $q_{2}(0)-q_{1}(0)=1$, $q_{3}(0)-q_{1}(0)=Q\,(0<Q<1)$ and $q_{2}(0)-q_{3}(0)=1-Q$, if the three bodies form a Eulerian collinear central configuration, then the following equality
\begin{eqnarray}\label{3.1}
\frac{m_{3}Q^{-2}+m_2}{m_{3}Q+m_2}=\frac{m_{1}+m_{3}(1-Q)^{-2}}{m_{1}+m_{3}(1-Q)}
\end{eqnarray}
holds.
\end{lemma}

Let us give another useful lemma.
\begin{lemma}\label{lem3.2}\cite[Lemma 2]{Perko1985}
Let $q_j(0)$ be given in Section 2, then
\begin{eqnarray}\label{3.2}
\sum_{j=1}^{N-1}\frac{1-q_j(0)}{|1-q_j(0)|^3}=\frac{1}{4}\sum_{j=1}^{N-1}csc(\frac{\pi j}{N}).
\end{eqnarray}
\end{lemma}

We also need some properties about eigenvalues $\lambda_{k}$ and eigenvectors $\nu_{k}$ for circulant matrices, and firstly we introduce the concept of circulant matrices.
\begin{definition}\cite[Page 65, lines 1-2 and Page 66, lines 16-17]{Marcus1964}
An $N\times N$ matrix $C=(c_{kj})$ is circulant if $c_{k,j}=c_{k-1,j-1}$, where $c_{0,j}$ and $c_{k,0}$ are identified with $c_{N,j}$ and $c_{k,N}$, respectively.
\end{definition}

Then we have
\begin{lemma}\label{lem3.3}\cite[Page 303, lines 2-7]{Perko1985}
Every circulant matrix $C$ has the same forms of the eigenvalues $\lambda_{k}(C)$ and the corresponding eigenvectors $\nu_{k}$, more precisely,
\begin{eqnarray}\label{3.3}
\lambda_{k}(C)=\sum_{j=1}^{N}c_{1,j}q_{k-1}^{j-1}(0), \,\,\,\,\nu_{k}=(q_{k-1}(0), q_{k-1}^{2}(0),\ldots,q_{k-1}^{N}(0))^{T}, \,\,\,k=1,2,\ldots,N,
\end{eqnarray}
where $q_{k-1}(0)=e^{i\theta_{k-1}}=e^{i\frac{2(k-1)\pi}{N}}$.
\end{lemma}

Let the matrix $A$ be as follows:
\begin{eqnarray}\label{3.4}
A=(a_{k,\,j}), \,\, a_{k,\,j}=
\left\{\begin{array}{ll} \frac{1-q_{j-k}(0)}{|1-q_{j-k}(0)|^{3}},& k\neq j, \\
\quad 0, & k=j.
\end{array}\right.
\end{eqnarray}
Obviously, $A$ is a circulant matrix. Then by Lemma \ref{lem3.3}, we have the following lemmas:

\begin{lemma}\label{lem3.4}\cite[Lemma 12]{Perko1985}
The eigenvalues of $A$ have the property that for $k\neq N$ and $N\geq4$, we have $\lambda_{k}\neq0$ except that $\lambda_{\frac{N+1}{2}}=0$ for odd $N$.
\end{lemma}

\begin{lemma}\label{lem3.5}\cite[Page 305, lines 13-14]{Perko1985}
The eigenvectors $\nu_{k} \, (k=1,2,\ldots,N)$ of $N\times N$ circulant matrix forms a basis of $\mathbb{C}^{N}$.
\end{lemma}

\begin{lemma}\label{lem3.6}\cite[Page 65, equality (4) and Page 66, equality (5) ]{Marcus1964}
If the conjugate transpose of $\nu_{k}$ is denoted by $\bar{\nu}_{k}^{T}$, then
\begin{eqnarray*}
\bar{\nu}_{k}^{T}\nu_{j}=\cases{N, \quad \   k=j, \cr 0, \quad \  \  k\neq j;}  \quad \  \,\,\,\,(q_{-1}(0),q_{-2}(0),\ldots,q_{-N}(0))^{T}\bar{\nu}_{N}=N.
\end{eqnarray*}
\end{lemma}

\begin{lemma}\label{lem3.7}\cite[Lemma 2.1]{Yu2012}
$\frac{1}{N}\sum_{1\leq j<k\leq N}\frac{1}{|q_{j}(0)-q_{k}(0)|}=\sum_{1\leq j< N}\frac{1-q_{j}(0)}{|1-q_{j}(0)|^{3}}$.
\end{lemma}

\section{The proof of main result} \label{sec4}\setcounter{equation}{0}
We choose the geometric center of the regular $N$-polygon as the origin of our coordinate system, and the proof of Theorem \ref{thm2.1}
is divided into two steps.

\textbf{Step I.} We prove that if the ($N$+1)-th particle locates at the geometric center of the regular $N$-polygon, then
all the masses at the vertices of the regular $N$-polygon must be equal to each other, and in the following, we divide the proof of \textbf{Step I} into three parts.

\textbf{Part 1: The case of $N=2$}.
It is well-known that for 3-body problem, there are only two kinds of central configurations: Eulerian collinear central configuration and Lagrangian equilateral-triangle central configuration. Combining the third body locates at the the geometric center of the regular $N$-polygon ($N=2$), then we have $Q=1/2$. Then by Lemma \ref{lem3.1}, we have
\begin{eqnarray*}
\frac{4m_{3}+m_{2}}{\frac{1}{2}m_{3}+m_{2}}=\frac{m_{1}+4m_{3}}{m_{1}+\frac{1}{2}m_{3}},
\end{eqnarray*}
which implies
\begin{eqnarray*}
m_{3}(4m_{1}+\frac{1}{2}m_{2})=m_{3}(4m_{2}+\frac{1}{2}m_{1}).
\end{eqnarray*}
Thus $m_{1}=m_{2}$. Note that in (\ref{1.2}), $q_{3}=0$ and $c_{0}=0$ for $N$\,=\,2, then from $m_{1}=m_{2}$, we have
\begin{eqnarray*}
m_{1}=m_{2}=4(\omega^{2}-m_{3})=\frac{\omega^{2}-m_{3}}{\frac{1}{4}\csc\frac{\pi}{2}},
\end{eqnarray*}
i.e. we complete the proof of \textbf{Step I} for $N=2$.

\textbf{Part 2: The case of $N=3$}.

The ($N$+1)-th body locates at the geometric center of the regular $N$-polygon, i.e.
$x_{4}=0$, so
\begin{eqnarray*}
c_{0}=\frac{\sum\limits_{j=1}^{4}m_{j}x_{j}(0)}{\sum\limits_{j=1}^{4}m_{j}}
=\frac{m_{4}x_{4}(0)+\sum\limits_{j=1}^{3}m_{j}q_{j}(0)}{\sum\limits_{j=1}^{4}m_{j}}
=\frac{\sum\limits_{j=1}^{3}m_{j}q_{j}(0)}{\sum\limits_{j=1}^{4}m_{j}}.
\end{eqnarray*}
Observing that
system (\ref{1.2}) is equivalent to
\begin{eqnarray}\label{4.1}
\left\{\begin{array}{ll}
\sum\limits_{{j\neq k\atop{{1\leq j\leq 3}}}}\frac{m_{j}m_{k}(q_{j}(0)-q_{k}(0))}{|q_{j}(0)-q_{k}(0)|^{3}}+\frac{m_{4}m_{k}(x_{4}(0)-q_{k}(0))}{|x_{4}(0)-q_{k}(0)|^{3}}
=-\omega^{2}m_{k}(q_{k}(0)-c_{0}), \quad \  \  k=1,2,3, \\ \sum\limits_{1\leq j\leq 3}\frac{m_{j}m_{4}}{|q_{j}(0)-x_{4}(0)|^{3}}(q_{j}(0)-x_{4}(0))=-\omega^{2} m_{4}(q_{4}(0)-c_{0}), \quad \  \  k=4.
\end{array}\right.
\end{eqnarray}
Observing that $x_{4}(0)=0$, $q_3(0)=1$ and $q_{-j}(0)=q_{-j+3k}(0)\,(j=1,2,3,\,k\in \mathbb{Z})$,  from (\ref{4.1}),  we obtain
\begin{eqnarray}\label{4.2}
\left(
     \begin{array}{cccc}
  0 & \frac{1-q_1(0)}{|1-q_1(0)|^{3}} & \frac{1-q_2(0)}{|1-q_2(0)|^{3}} &  1  \\
 \frac{1-q_2(0)}{|1-q_2(0)|^{3}} & 0 & \frac{1-q_1(0)}{|1-q_1(0)|^{3}} &  1 \\
  \frac{1-q_1(0)}{|1-q_1(0)|^{3}} & \frac{1-q_2(0)}{|1-q_2(0)|^{3}} & 0 &  1  \\
  q_{1}(0) & q_{2}(0) & q_{3}(0) & 0  \\
 \end{array}
   \right)*\left(
     \begin{array}{cccccc}
  m_1  \\
  m_2  \\
 m_3 \\
m_4 \\
 \end{array}
   \right)=\left(
     \begin{array}{cccccc}
  \omega^{2}-\omega^{2} c_{0}q_{2}(0)  \\
  \omega^{2}-\omega^{2} c_{0}q_{1}(0)  \\
 \omega^{2}-\omega^{2} c_{0} \\
\omega^{2} c_{0} \\
 \end{array}
   \right)
\end{eqnarray}
Therefore
\begin{eqnarray}\label{4.3}
&&\Big[\frac{1-q_1(0)}{|1-q_1(0)|^{3}}m_{1}+\frac{1-q_2(0)}{|1-q_2(0)|^{3}}m_{2}+m_{4}\Big]+(m_1q_1(0)+m_2q_2(0)+m_3q_3(0))
\nonumber \\ &=&[\omega^{2}-\omega^{2} c_{0}]+\omega^{2} c_{0}=\omega^{2}\in\mathbb{R}.
\end{eqnarray}
Since $q_1(0)=e^{i\frac{2\pi}{3}}$ and $q_2(0)=e^{i\frac{4\pi}{3}}$, then $Re(q_1(0))=Re(q_2(0))$ and $Im(q_1(0))=-Im(q_2(0))=\sqrt{3}i/2$. Therefore, $1/|1-q_1(0)|^{3}=1/|1-q_2(0)|^{3}$.
Note that $m_{1},m_{2},m_{3},m_{4}\in\mathbb{R}$, then from (\ref{4.3}),
\begin{eqnarray*}
-\frac{\sqrt{3}}{2}i\Big[\frac{1}{|1-q_1(0)|^{3}}(m_1-m_2)\Big] +\frac{\sqrt{3}}{2}i(m_1-m_2)=0,
\end{eqnarray*}
which implies that $m_1=m_2=m$.

Next we will prove that $m_3=m$. In fact, by (\ref{4.2}), we have
\begin{eqnarray}\label{4.4}
\left\{\begin{array}{ll}
(m_3-m)\frac{1-q_2(0)}{|1-q_2(0)|^{3}}+(m-m_3)\frac{1-q_1(0)}{|1-q_1(0)|^{3}}=\omega^{2} c_{0}(q_1(0)-q_2(0)), \\ mq_1(0)+mq_2(0)+m_3=\omega^{2} c_0.
\end{array}\right.
\end{eqnarray}
From the representations of $q_1(0)$ and $q_2(0)$, then $q_1(0)^2=q_2(0)$ and $q_2(0)^2=q_1(0)$. It follows from the second identity in (\ref{4.4}) that
\begin{eqnarray*}
\omega^{2} c_{0}(q_1(0)-q_2(0)) &=& [mq_1(0)+mq_2(0)+m_3](q_1(0)-q_2(0))
\\ &=& m(q_1(0)^{2}-q_2(0)^{2})+m_3(q_1(0)-q_2(0))
\\ &=&(m_{3}-m)(q_1(0)-q_2(0)).
\end{eqnarray*}
 Then combining the first part of system (\ref{4.4}), we obtain
\begin{eqnarray*}
&&(m_3-m)\frac{1-q_2(0)}{|1-q_2(0)|^{3}}+(m-m_3)\frac{1-q_1(0)}{|1-q_1(0)|^{3}}-(m_3-m)(q_1(0)-q_2(0))
\\ &=&(m_3-m)\Big[\frac{1}{|1-q_1(0)|^{3}}-1\Big](q_1(0)-q_2(0))
\\ &=&(m_3-m)\Big[\frac{1}{3\sqrt{3}}-1\Big]\sqrt{3}i=0,
\end{eqnarray*}
which suggests that $m_3=m$.

Note that $x_{4}=c_0=0$ ($N$=3) and $m_{1}=m_{2}=m_{3}$, by choosing $k=3$ in (\ref{1.2}), we have
\begin{eqnarray*}
m\sum_{j=1}^{2}\frac{(q_j(0)-1)}{|q_j(0)-1|^{3}}-m_4
=-\omega^{2},
\end{eqnarray*}
and thus
\begin{eqnarray*}
m_{1}=m_{2}=m_{3}=\frac{\omega^{2}-m_{4}}{\frac{1}{4}\sum_{j=1}^{2}\csc\frac{j\pi}{3}},
\end{eqnarray*}
if one uses the identity (\ref{3.2}). Hence we complete the proof of \textbf{Step I} for $N=3$.

\textbf{Part 3: The case of $N\geq4$}.

By $x_{N+1}(0)=0$ and $|q_{j}(0)|=|e^{i\theta_{j}}|=1$ with $j=1,2,\ldots,N$, we gain from system (\ref{1.2}) that
\begin{eqnarray}\label{4.5}
\left\{\begin{array}{ll}
\sum\limits_{{j\neq k\atop{{1\leq j\leq N}}}}\frac{m_{j}}{|q_{j}(0)-q_{k}(0)|^{3}}(q_{j}(0)-q_{k}(0))-m_{N+1}q_{k}(0)=-\omega^{2} \Big[q_{k}(0)-\frac{\sum\limits_{j=1}^{N}m_{j}q_{j}(0)}{\sum\limits_{j=1}^{N+1}m_{j}}\Big], \ \ k\neq N+1, \\ \sum\limits_{1\leq j\leq N}m_{j}q_{j}(0)=\omega^{2} \frac{\sum\limits_{j=1}^{N}m_{j}q_{j}(0)}{\sum\limits_{j=1}^{N+1}m_{j}},\ \  k= N+1.
\end{array}\right.
\end{eqnarray}
From the second equation in (\ref{4.5}), then $\sum\limits_{j=1}^{N}m_{j}q_{j}(0)=0$ or $\omega^{2}=\sum\limits_{j=1}^{N+1}m_{j}$.
In the following, we will discuss the two cases.

\textbf{Case 2.1. $\omega^{2}=\sum\limits_{j=1}^{N+1}m_{j}$}.

By the first equation in (\ref{4.5}), we have
\begin{eqnarray*}
\sum\limits_{{j\neq k\atop{{1\leq j\leq N}}}}\frac{q_{j}(0)-q_k(0)}{|q_{j}(0)-q_{k}(0)|^{3}}m_{j}=-\sum_{j=1}^{N}m_{j}q_{k}(0)+\sum_{j=1}^{N}m_{j}q_{j}(0), \quad \  \  k\neq N+1,
\end{eqnarray*}
which suggests that
\begin{eqnarray}\label{4.6}
\sum\limits_{{j\neq k\atop{{1\leq j\leq N}}}}\frac{1-q_{j-k}(0)}{|1-q_{j-k}(0)|^{3}}m_{j}=\sum_{j=1}^{N}m_{j}-
\sum_{j=1}^{N}m_{j}q_{j}(0)q_{-k}(0), \quad \  \  k\neq N+1.
\end{eqnarray}

Let the circulant matrix $A$ be given by (\ref{3.4}). With the help of (\ref{4.6}) and Lemma \ref{lem3.3}, then
\begin{eqnarray}\label{4.7}
A(m_{1},m_{2},\ldots,m_{N})^{T}=(\sum_{j=1}^{N}m_{j})\nu_{1}-
(\sum_{j=1}^{N}m_{j}q_{j}(0))\nu,
\end{eqnarray}
where $\nu=(q_{-1}(0),q_{-2}(0),\ldots,q_{-N}(0))$.

Observing that $q_{j}(0)=e^{i\theta_{j}}=e^{i\frac{2j\pi}{N}}$, it yields that $q_{-j}(0)=q_{N-1}^j(0)$. Therefore (\ref{4.7}) is equivalent to
\begin{eqnarray}\label{4.8}
A(m_{1},m_{2},\ldots,m_{N})^{T}=(\sum_{j=1}^{N}m_{j})\nu_{1}-
(\sum_{j=1}^{N}m_{j}q_{j}(0))\nu_N,
\end{eqnarray}
where $\nu_1$ and $\nu_N$ are defined in (\ref{3.3}).

On the other hand, from Lemma \ref{lem3.5}, there exist $\alpha_{1},\alpha_{2},\ldots,\alpha_{N}\in\mathbb{C}$ such that
\begin{eqnarray}\label{4.9}
(m_{1},m_{2},\ldots,m_{N})^{T}=\alpha_{1}\nu_{1}+\alpha_{2}\nu_{2}+\ldots+\alpha_{N}\nu_{N},
\end{eqnarray}
where $\nu_{1},\nu_{2},\ldots,\nu_{N}$ are defined in (\ref{3.1}).

Employing (\ref{4.8}) and (\ref{4.9}), there exist eigenvalues $\lambda_{1},\lambda_{2},\ldots,\lambda_{N}\in\mathbb{C}^{N}$ such that
\begin{eqnarray}\label{4.10}
A(\alpha_{1}\nu_{1}+\alpha_{2}\nu_{2}+\ldots+\alpha_{N}\nu_{N})
&=&\alpha_{1}A\nu_{1}+A\alpha_{2}\nu_{2}+\ldots+A\alpha_{N}\nu_{N}
\nonumber\\
&=&\lambda_{1}\alpha_{1}\nu_{1}
+\lambda_{2}\alpha_{2}\nu_{2}+\ldots+\lambda_{N}\alpha_{N}\nu_{N}\nonumber\\
&=&(\sum_{j=1}^{N}m_{j})\nu_{1}-(\sum_{j=1}^{N}m_{j}q_{j}(0))\nu_{N}.
\end{eqnarray}
From Lemma \ref{lem3.6}, we have $\nu_{1},\nu_{2},\ldots,\nu_{N}$ is linearly independent. Then combining (\ref{4.10}), we see that $\lambda_{2}\alpha_{2}=\lambda_{3}\alpha_{3}=\ldots=\lambda_{N-1}\alpha_{N-1}=0$. In fact, by Lemma \ref{lem3.4}, we have $\lambda_{k}\neq0$ except that $\lambda_{\frac{N+1}{2}}=0$, which implies
\begin{eqnarray}\label{4.11}
\alpha_{k}=0, \quad\ \  where \,\, k=2,3,\ldots,N-1 \,\, and \,\, k\neq\frac{N+1}{2}.
\end{eqnarray}
Then we have two subcases.

\textbf{Case 2.1.1. $N$ is a even number}.

 Combining (\ref{4.9}) with (\ref{4.11}), we have
\begin{eqnarray}\label{4.12}
(m_{1},m_{2},\ldots,m_{N})^{T}-\alpha_{1}\nu_{1}=\alpha_{N}\nu_{N},
\end{eqnarray}
where $\nu_{1}=(1,1,\ldots,1)^{T}$ and $\nu_{N}=(q_{N-1}(0),q_{N-2}(0),\ldots,q_{0}(0))^{T}$.
We set $\alpha_1=a_1+ib_1$ and $\alpha_N=a_2+ib_2$. Notice that $(m_{1},m_{2},\ldots,m_{N})^{T}$ is a real vector. From (\ref{4.12}), then $(a_2+ib_2)q_j+ib_1$ are real numbers for $ 0\leq j\leq N-1$. Observing that
\begin{eqnarray*}
(a_2+ib_2)q_j+ib_1&=&(a_2+ib_2)e^{i\frac{2j\pi}{N}}+ib_1
\\ &=& a_2\cos(\frac{2j\pi}{N})-b_2\sin(\frac{2j\pi}{N})+
i\Big[b_2\cos(\frac{2j\pi}{N})+a_2\sin(\frac{2j\pi}{N})+b_1\Big],
\end{eqnarray*}
therefore
\begin{eqnarray*}
b_2\cos(\frac{2j\pi}{N})+a_2\sin(\frac{2j\pi}{N})+b_1=0,\  \ j=0,1,_{\cdots},N-1.
\end{eqnarray*}
By choosing $j=0$ and $j=N/2$, we conclude that
\begin{eqnarray*}
b_2+b_1=0,\  \  -b_2+b_1=0.
\end{eqnarray*}
Thus $b_2=b_1=0$. So from (\ref{4.12}), it yields that
\begin{eqnarray*}
(m_{1}-a_1,m_2-a_1,\ldots,m_N-a_1)^{T}=a_2\nu_{N},
\end{eqnarray*}
which suggests that $a_2e^{i\frac{2j\pi}{N}}$ are real numbers for $j=0,1,\ldots,N-1$. Since $N\geq 4$, then $a_2=0$.
Hence, $(m_{1},m_{2},\ldots,m_{N})^{T}=a_1\nu_{1}$, i.e., $m_{1}=m_{2}=\ldots=m_{N}$.

\textbf{Case 2.1.2. $N$ is a odd number}.

By (\ref{4.9}) and (\ref{4.11}), we have
\begin{eqnarray*}
(m_{1},m_{2},\ldots,m_{N})^{T}=\alpha_{1}\nu_{1}+\alpha_{N}\nu_{N}+\alpha_{\frac{N+1}{2}}\nu_{\frac{N+1}{2}}.
\end{eqnarray*}
So \begin{eqnarray}\label{4.13}
(m_{1},m_{2},\ldots,m_{N})^{T}-\alpha_{1}\nu_{1}=\alpha_{N}\nu_{N}+\alpha_{\frac{N+1}{2}}\nu_{\frac{N+1}{2}}.
\end{eqnarray}
where $\nu_{1}=(1,1,\ldots,1)^{T}$, $\nu_{N}=(q_{N-1}(0),q_{N-2}(0),\ldots,q_{0}(0))^{T}$ and
\begin{eqnarray*}
\nu_{\frac{N+1}{2}}=(q_{\frac{N+1}{2}-1}(0),q_{\frac{N+1}{2}-1}^2(0),\ldots,
q_{\frac{N+1}{2}-1}^N(0))^{T}.
\end{eqnarray*}
 We set $\alpha_1=a_1+ib_1$, $\alpha_N=a_2+ib_2$ and  $\alpha_{\frac{N+1}{2}}=a_3+ib_3$.
By (\ref{4.13}),  $(a_2+ib_2)q_{N-j}(0)+(a_3+ib_3)q_{\frac{N+1}{2}-1}^j(0)+ib_1$ are real numbers for $1\leq j\leq N$. Notice that
\begin{eqnarray*}
&&(a_2+ib_2)q_{N-j}(0)+(a_3+ib_3)q_{\frac{N+1}{2}-1}^j(0)+ib_1
\\ &=&(a_2+ib_2)e^{i\frac{2(N-j)\pi}{N}}+(a_3+ib_3)e^{i\frac{(N-1)j\pi}{N}}
+ib_1
\\ &=& a_2\cos(\frac{2j\pi}{N})+b_2\sin(\frac{2j\pi}{N})+(-1)^j(a_3\cos(\frac{k\pi}{N})+b_3\sin(\frac{j\pi}{N}))
\\&&+i\Big[b_2\cos(\frac{2j\pi}{N})-a_2\sin(\frac{2j\pi}{N})+(-1)^j(b_3\cos(\frac{j\pi}{N})-a_3\sin(\frac{j\pi}{N}))+b_1\Big],
\end{eqnarray*}
therefore
\begin{eqnarray*}
b_2\cos(\frac{2j\pi}{N})-a_2\sin(\frac{2j\pi}{N})+(-1)^j(b_3\cos(\frac{j\pi}{N})-a_3\sin(\frac{j\pi}{N}))+b_1=0,\  \ j=1,2,\ldots,N.
\end{eqnarray*}
Observing that $N\geq 4$ and $N$ is odd, so $N\geq 5$. By choosing $j=N, N-1, N-2$, $2$ and $1$, we conclude the following linear equations
\begin{eqnarray}\label{4.14}
\cases{b_2+b_3+b_1=0, \cr b_2\cos(\frac{2\pi}{N})+a_2\sin(\frac{2\pi}{N})-b_3\cos(\frac{\pi}{N})
-a_3\sin(\frac{\pi}{N})+b_1=0,
\cr b_2\cos(\frac{4\pi}{N})+a_2\sin(\frac{4\pi}{N})+
b_3\cos(\frac{2\pi}{N})
+a_3\sin(\frac{2\pi}{N})+b_1=0,
\cr b_2\cos(\frac{4\pi}{N})-a_2\sin(\frac{4\pi}{N})+b_3\cos(\frac{2\pi}{N})
-a_3\sin(\frac{2\pi}{N})+b_1=0,
\cr
b_2\cos(\frac{2\pi}{N})-a_2\sin(\frac{2\pi}{N})-
b_3\cos(\frac{\pi}{N})+a_3\sin(\frac{\pi}{N})+b_1=0.}
\end{eqnarray}
We solve from (\ref{4.14}) that
\begin{eqnarray*}
b_1=b_2=b_3=a_2=a_3=0.
\end{eqnarray*}
Therefore from (\ref{4.13}), $\alpha_N=a_2+ib_2$ and $\alpha_{\frac{N+1}{2}}=a_3+ib_3$, we arrive at
$(m_{1},m_{2},\ldots,m_{N})^{T}=a_1\nu_{1}$, i.e. $m_{1}=m_{2}=\ldots=m_{N}$.

Thus if $\omega^{2}=\sum\limits_{j=1}^{N+1}m_{j}$, we obtain that $m_{1}=m_{2}=\ldots=m_{N}$.

\textbf{Case 2.2. $\sum\limits_{j=1}^{N}m_{j}q_{j}(0)=0$}.

Remind the first part of system (\ref{4.5}), we have
\begin{eqnarray*}
\sum\limits_{{j\neq k\atop{{1\leq j\leq N}}}}\frac{q_{j}(0)-q_{k}(0)}{|q_{j}(0)-q_{k}(0)|^{3}}m_{j}=(-\omega^{2}+m_{N+1})q_{k}(0), \quad \  \  k\neq N+1,
\end{eqnarray*}
which implies
\begin{eqnarray*}
\sum\limits_{{j\neq k\atop{{1\leq j\leq N}}}}\frac{1-q_{j-k}(0)}{|1-q_{j-k}(0)|^{3}}m_{j}=\omega^{2}-m_{N+1}, \quad \  \  k\neq N+1.
\end{eqnarray*}
Therefore we know that
\begin{eqnarray}\label{4.15}
A(m_{1},m_{2},\ldots,m_{N})^{T}=(\omega^{2}-m_{N+1})\nu_{1}.
\end{eqnarray}
Employing (\ref{4.9}) and (\ref{4.15}), there exist eigenvalues $\lambda_{1},\lambda_{2},\ldots,\lambda_{N}\in\mathbb{C}$ such that
\begin{eqnarray*}
A(\alpha_{1}\nu_{1}+\alpha_{2}\nu_{2}+\ldots+\alpha_{N}\nu_{N})
&=&\alpha_{1}A\nu_{1}+A\alpha_{2}\nu_{2}+\ldots+A\alpha_{N}\nu_{N}
\nonumber\\
&=&\lambda_{1}\alpha_{1}\nu_{1}
+\lambda_{2}\alpha_{2}\nu_{2}+\ldots+\lambda_{N}\alpha_{N}\nu_{N}\nonumber\\
&=&(\omega^{2}-m_{N+1})\nu_{1}.
\end{eqnarray*}
The left process are similar to \textbf{Case 2.1}, then we also obtain $m_{1}=m_{2}=\ldots=m_{N}$.

From \textbf{Cases 2.1} and \textbf{2.2}, we obtain that $m_{1}=m_{2}=\ldots=m_{N}$, and by our choice for the origin of coordinate system, we know the mass center of the regular
$N$-polygon coincides with the geometric center, furthermore, $x_{N+1}(0)=0$, $\sum_{j=1}^{N}m_{j}q_{j}(0)/\sum_{j=1}^{N}m_{j}=0$,
therefore $c_{0}=0$. Then combining (\ref{1.2}), $x_{N+1}=0$, $c_{0}=0$, $m_{1}=m_{2}=\ldots=m_{N}$ and Lemma \ref{lem3.2}, we have
\begin{eqnarray*}
m_{1}=m_{2}=\ldots=m_{N}=\frac{\omega^{2}-m_{N+1}}{\frac{1}{4}\sum_{j=1}^{N-1}\csc\frac{j\pi}{N}}.
\end{eqnarray*}
Therefore we complete the proof of \textbf{Step I} for $N\geq4$.

\textbf{Step II.} Under the assumption that the ($N$+1)-th ($N\geq2$) body locates at the geometric center of the regular $N$-polygon, we prepare to prove that if
\begin{eqnarray*}
m_{1}=m_{2}=\ldots=m_{N}=\frac{\omega^{2}-m_{N+1}}{\frac{1}{4}\sum_{j=1}^{N-1}\csc\frac{j\pi}{N}},
\end{eqnarray*}
then all the $N$+1 bodies can form a central configuration.

In fact, from system (\ref{1.2}), we see that under the assumption that the ($N$+1)-th body locates at the geometric center of the regular $N$-polygon, if (\ref{4.5}) holds for $N\geq2$, then all the $N$+1 bodies can form a central configuration.

We notice that for the configuration of planar ($N$+1)-body ($N\geq2$) problem with a regular $N$-polygon, if
\begin{eqnarray*}
m_{1}=m_{2}=\ldots=m_{N}=\frac{\omega^{2}-m_{N+1}}{\frac{1}{4}\sum_{j=1}^{N-1}\csc\frac{j\pi}{N}}=:m,
\end{eqnarray*}
and the ($N$+1)-th body locates at the geometric center of the regular $N$-polygon, then the mass center $c_{0}=0$, which implies the second part of (\ref{4.5}) always hold for $N\geq2$. So in the following, it suffice to prove that the first part of (\ref{4.5}) always hold for $N\geq2$.

There is no loss of generality in assuming $m_{N+1}=bm$ $(b>0)$, then combining $c_{0}=0$, $x_{i}(0)=q_{i}(0)$ $(i=1,2,\ldots,N)$, $x_{N+1}(0)=0$, \textbf{Definition 1.2} and the equivalent definitions of central configuration in \cite[Page 109, lines 1-7]{Llibre2015} imply that
\begin{eqnarray*}
\omega^{2}=\frac{U(\vec{x}(0))}{I(\vec{x}(0))},
\end{eqnarray*}
where
\begin{eqnarray*}
U(\vec{x}(0))=\sum_{1\leq k<j\leq N+1}\frac{m_j m_k}{|x_j(0)-x_k(0)|},
 \,\,\,\,I(\vec{x}(0))=\sum_{1\leq j\leq N+1}m_j|x_j(0)-c_{0}|^{2},
\end{eqnarray*}
we have
\begin{eqnarray*}
\omega^{2}&=&\Big[\sum_{1\leq k<j\leq N+1}\frac{m_j m_k}{|x_j(0)-x_k(0)|}\Big]\times\Big[\frac{1}{\sum\limits_{1\leq j\leq N+1}m_j|x_j(0)-c_{0}|^{2}}\Big]
\nonumber\\ &=&\frac{1}{(\sum\limits_{1\leq j\leq N}m|q_j(0)|^{2})}\Big(\sum_{1\leq k<j\leq N}\frac{m^{2}}{|q_j(0)-q_k(0)|}+\sum_{1\leq k\leq N}\frac{bm^{2}}{|q_k(0)|}\Big)
\nonumber\\&=&\frac{1}{N}\Big(\sum_{1\leq k<j\leq N}\frac{m}{|q_j(0)-q_k(0)|}+Nbm\Big)
\nonumber\\&=&\frac{m}{N}\sum_{1\leq k<j\leq N}\frac{1}{|q_j(0)-q_k(0)|}+bm.
\end{eqnarray*}
With the help of Lemma \ref{lem3.7}, then
\begin{eqnarray}\label{4.16}
\omega^{2}=\frac{m}{N}\sum_{1\leq k<j\leq N}\frac{1}{|q_j(0)-q_k(0)|}+bm=m\sum_{1\leq j< N}\frac{1-q_{j}(0)}{|1-q_{j}(0)|^{3}}+bm.
\end{eqnarray}
Combining (\ref{4.16}) and the fact $m_{1}=m_{2}=\ldots=m_{N}=m_{N+1}/b=m$, it can be checked that the first equations in (\ref{4.5}) always hold for $N\geq2$. Thus all the $N$+1 bodies can form a central configuration.

By now, from \textbf{Step I} and \textbf{Step II} , we complete the proof of Theorem \ref{thm2.1}.  $\Box$

\vskip5mm\noindent
\textbf{Acknowledgements.}
\vskip2mm\par
The authors would like to thank Professor Zhifu Xie for some discussions. The first author is partially supported by research funding project of Guizhou Minzu University (GZMU[2019]QN04). The second author is partially supported by NSF of China (11501577). The third author is partially supported by NSF of China (116712787).

\end{document}